\documentclass[11pt,notitlepage,a4paper,reqno]{amsart}
\usepackage{amsmath}
\usepackage{amsfonts}
\usepackage{amssymb}
\usepackage{epsfig}
\usepackage[usenames]{color}
\usepackage{colortbl}
\newcommand{\RR}{\ensuremath{\mathbb{R}}}
\newcommand{\rn}{\ensuremath{\mathbb{R}^n}}
\newcommand{\NN}{\ensuremath{\mathbb{N}}}

\newcommand{\beq}{\begin{equation}}
\newcommand{\eeq}{\end{equation}}

\newcommand{\medint}{-\kern  -,375cm\int}

\def\x{\bar x}
\theoremstyle{plain}
\newtheorem{theorem}{Theorem}[section]
\newtheorem{corollary}[theorem]{Corollary}
\newtheorem{lemma}[theorem]{Lemma}
\newtheorem{proposition}[theorem]{Proposition}
\theoremstyle{plain}
\newtheorem{definition}[theorem]{Definition}
\theoremstyle{remark}
\newtheorem{remark}[theorem]{Remark}
\theoremstyle{plain}

\theoremstyle{plain}

\theoremstyle{plain}
\newtheorem*{cp*}{C.P}

\numberwithin{equation}{section}

\makeatletter

\renewcommand{\p@enumi}{\thesection.}
\makeatother

\allowdisplaybreaks

\begin{document}

\baselineskip.5cm


\title[Mean width rearrangements]{Combination and mean width rearrangements of solutions of elliptic
equations in convex sets}
\author[Paolo Salani]{Paolo Salani}
\address{Dipartimento di Matematica ``U. Dini'' -
Universit\`a di Firenze - viale Morgagni 67/A, 50134 Firenze
(Italy)} \email{paolo.salani@unifi.it}

\begin{abstract}
We introduce a method to compare solutions of different equations in different domains.
As a consequence, we define a new kind of rearrangement which applies to solution of fully nonlinear equations
$F(x,u,Du,D^2u)=0$, not necessarily in divergence form, in convex domains and we obtain Talenti's type results for this kind of rearrangement. 
\end{abstract}

\maketitle


\section{Introduction}
Rearrangements are among the most powerful tools in analysis.
Roughly speaking they manipulate the shape of an object while preserving someone of its relevant geometric properties. Typically, a rearrangement of a function is performed by acting separately on each of its level sets. Probably the most famous one is the radially symmetric decreasing rearrangement, or {\em Schwarz symmetrization}: the {\em Schwarz symmetrand} of a continuous function $w\geq 0$ is the function $w^\star$ whose superlevel sets are concentric balls (usually centered at the origin) with the same measure of the corresponding superlevel sets of $w$.  Notice that $w^\star$, by definition, is equidistributed with $w$. When applied to the study of solutions of partial differential equations with a divergence structure, this usually leads to a comparison between the solution  in a generic domain and the solution of (a possibly "rearranged" version of) the same equation in a ball with the same measure of the original domain. An archetypal result of this type is the following (see \cite{Talenti}): let $u^\star$ be the Schwarz symmetrand of the solution $u$ of
\begin{equation}\label{Talentieq}
\left\{\begin{array}{ll}
\Delta u+f(x)= 0\quad&
\mbox{in }\Omega\\
\\
u=0\quad&\mbox{on }\partial\Omega
\end{array}
\right.
\end{equation}
and let $v$ be the solution of
$$
\left\{\begin{array}{ll}
\Delta v+f^\star(x)=0\quad&
\mbox{in }\Omega^\star\,,\\
\\
v=0\quad&\mbox{on }\partial\Omega^\star\,,
\end{array}
\right.
$$
where $\Omega^\star$ is the ball (centered at the origin) with the same measure of $\Omega$, $f$ is a non-negative function and $f^\star$ is the Schwarz symmetrand of $f$. Then, under suitable summability assumptions on $f$, it holds
\begin{equation}\label{comparisonTalenti}
u^\star\leq v\quad\mbox{in }\Omega^\star\,.
\end{equation}
whence
\begin{equation}\label{comparisonTalenti1}
\|u\|_{L^p(\Omega)}\leq\|v\|_{L^p(\Omega^\star)}
\end{equation}
for every $p>0$, including $p=+\infty$.

Actually the above Talenti's comparison principle (\ref{comparisonTalenti})-(\ref{comparisonTalenti1}) applies to more general situations and the Laplace operator in
(\ref{Talentieq}) can be substituted by operators like
$$
\mbox{div}(a_{ij}(x)u_j)+c(x)u
$$
or even more general ones (see for instance \cite{ALT, ATLM, Talenti, Talenti2, TV}), 
but always {\em  in divergence form}.

Here we introduce a new kind of rearrangement, which permits to obtain comparison results similar to (\ref{comparisonTalenti})-(\ref{comparisonTalenti1})  {\em for very general equations, not necessarily in divergence form}, between a classical solution in a convex domain $\Omega$ and the solution in the ball $\Omega^\sharp$ with the same mean width of $\Omega$. Recall that {\em the mean width} $w(\Omega)$ of $\Omega$ is defined as follows:
$$w(\Omega)=\frac{1}{n\omega_{n}} \int_{S^{n-1}} \! \big(h(\Omega,\xi)+h(\Omega,-\xi)\big) \, d\xi=\frac{2}{n\omega_{n}} \int_{S^{n-1}} \! h(\Omega,\xi) \, d\xi\,,
$$
where $h(\Omega, \cdot)$ is the support function of $\Omega$ (then $w(\Omega,\xi)=w(\Omega,-\xi)=h(\Omega,\xi)+h(\Omega,-\xi)$ is {\em the width } of $\Omega$ in direction $\xi$ or $-\xi$) and $\omega_{n}$ is the measure of the unit ball in $\mathbb{R}^{n}$. When $\Omega$ is a ball, $w(\Omega)$ simply coincides with its diameter; in the plane $w(\Omega)$ coincides with the perimeter of $\Omega$, up to a factor $\pi^{-1}$. See Section 2 for more details, notation and definitions. 

Precisely, we will deal with problems of the following type
\begin{equation}
\label{pb1}
\left\{
\begin{array}{ll}
F(x,u,D u, D^2u)=0 & \textrm{ in } \Omega\,,\\
u=0 & \textrm{ on } \partial \Omega\,,\\
u>0 &\textrm{ in }\Omega\,,
\end{array}
\right.
\end{equation}
where $F(x,t,\xi,A)$ is a continuous proper elliptic operator acting on $\RR^n\times
\RR\times\RR^n\times S_n$ and $\Omega$ is an open bounded convex subset of $\RR^n$.
Here $D u$ and $D^2u$ are the gradient and the Hessian matrix of the
function $u$ respectively, $S_n$ is the set of the $n\times n$
real symmetric matrices.

We will see how, given a solution ${u}$ of problem (\ref{pb1}) and a parameter $p>0$, it is possible to associate to ${u}$ a symmetrand $u_p^\sharp$ which is defined in a ball $\Omega^\sharp$ having the same mean width of $\Omega$ and, under suitable
assumptions on the operator $F$ (see Theorem \ref{thmrearr}), we obtain a pointwise comparison analogous to (\ref{comparisonTalenti}) between $u_p^\sharp$ and the solution ${v}$ in $\Omega^\sharp$, that is
\begin{equation}\label{comparisonSalani}
u_p^\sharp\leq {v}\quad\mbox{ in }\Omega^\sharp\,,
\end{equation}
where $v$ is the solution of
\begin{equation}
\label{pbsharp}
\left\{
\begin{array}{ll}
F(x,v,D v, D^2v)=0 & \textrm{ in } \Omega^\sharp\,,\\
v=0 & \textrm{ on } \partial \Omega^\sharp\,,\\
v>0 &\textrm{ in }\Omega^\sharp\,,
\end{array}
\right.
\end{equation}
Then from (\ref{comparisonSalani}) we get
\begin{equation}\label{lpustargequ}
\|{u}\|_{L^q(\Omega)}\leq \|{v}\|_{L^q(\Omega^\sharp)}\quad\mbox{for every }q\in(0,+\infty]\,.
\end{equation}
The  precise definition of $u_p^\sharp$ is actually quite involved and it will be given in Section 5. Here we just say that $u_p^\sharp$ is not equidistributed with ${u}$, in contrast with Schwarz symmetrization; indeed the measure of the super level sets of $u_p^\sharp$ is greater than the measure of the corresponding super level sets of ${u}$. 

The results of this paper are based on the refinement of a technique developed in \cite{BS, cuoghis, IS5} (and inspired by \cite{all}) to study concavity properties of solution of elliptic and parabolic equations in convex rings and in convex domains.
It is shown here that this refinement permits {\em to compare solutions of different equations in different domains} and this is in fact the main result of the paper, see Theorem \ref{mainthm}. More explicitly, consider two convex sets $\Omega_0$ and $\Omega_1$ and a real number $\mu\in(0,1)$, and denote by $\Omega_\mu$ the {\em Minkowski convex combination} (with coefficient $\mu$) of $\Omega_0$ and $\Omega_1$, that is
$$
\Omega_\mu=(1-\mu)\Omega_0+\mu\,\Omega_1=\{(1-\mu) x_0+\mu\,x_1\,:\, x_0\in\Omega_0,\, x_1\in\Omega_1\}\,.
$$
Correspondingly, let $u_0$, $u_1$ and $u_\mu$ be the solutions of 
$$
(P_i)\quad
\left\{
\begin{array}{ll}
F_i(x,u_i,D u_i, D^2u_i)=0 & \textrm{ in } \Omega_i\,,\\
u_i=0 & \textrm{ on } \partial \Omega_i\,,\qquad i=0,1,\mu\\
u_i>0 &\textrm{ in }\Omega_i\,.
\end{array}
\right.
$$
Roughly speaking (the precise statement will be given in Section 4)
Theorem \ref{mainthm} states that, under suitable assumptions on the operators $F_0$, $F_1$ and $F_\mu$, it is possible to compare $u_\mu$ with 
a suitable convolution of $u_0$ and $u_1$. 
Such a result has obviously its own interest and it has several interesting consequences, among which there is the rearrangements technique sketched above. 
\medskip

The paper is organized as follows. In Section 2 we introduce notation and recall some useful  notions and known results.
Section 3 is dedicated to the so-called $(p,\mu)$-convolution of non-negative functions. In Section 4 it is stated Theorem \ref{mainthm}, the main theorem of the paper, which is proved in Section 5. Section 6 is devoted to rearrangements: it contains the definition of $u_p^\sharp$ and Theorem \ref{thmrearr}. In Section 7 some examples and applications are presented.

\section{Notation and preliminaries}

For $A\subseteq\RR^n$, we denote by $\overline A$, $\partial A$ and $|A|$ its
closure, its boundary and its measure, respectively.

Let $n\geq 2$, $x\in \RR^n$
and $r>0$: $B(x,r)$ is the euclidean ball of radius $r$ centered
at $x$, i.e.
$$B(x,r)=\left\{z\in \RR^n\,:\,|z-x|<r\right\}.$$
In particular we set $B=B(0,1)$, $S^{n-1}=\partial B$ and $\omega_n=|B|$.

We denote by $S_n$ the space of $n\times n$ real symmetric matrices and by $S_n^+$ and $S_n^{++}$ the cones of nonnegative and positive definite symmetric matrices, respectively. If $A,B\in S_n$, by $A\geq 0\,(>0)$ we mean that $A\in S_n^+\,(S_n^{++})$ and $A\geq B$ means $A-B\geq 0$.

$SO(n)$ is the special orthogonal group of $\rn$, that is the space of rotations in $\rn$, i.e. $n\times n$ orthogonal matrices with determinant $1$. 


With the symbol $\otimes$
we denote the direct product between vectors in $\RR^n$, that is, for $x=(x_1,\dots,x_n)$ and $y=(y_1,\dots,y_n)$,
$x\otimes y$ is the $n\times n$ matrix  with entries $(x_iy_j)$ for
$i,j=1,...,n$.

\subsection{Viscosity solutions}

We will make use of basic viscosity techniques; here we recall only few notions and we refer to the User's Guide
\cite{UG} and to the books \cite{caffacabre, Koike} for more details.

The continuous operator $F:\RR^n\times\RR\times \RR^n\times S_n\rightarrow \RR$
is said \emph{proper} if
$$
F(x,r,\xi ,A)\leq F(x,s,\xi ,A) \textrm{ whenever }r\geq s.
$$
Let $\Gamma$ be a convex cone in $S_n$, with vertex at the origin and
containing the cone of nonnegative definite symmetric matrices $S^+_n$.
We say that $F$ is \emph{degenerate elliptic} in $\Gamma$ if
$$
F(x,u,\xi ,A)\leq F(x,u,\xi ,B) \textrm{ whenever } A\leq B,\, A,B\in\Gamma\,.
$$

We set $\Gamma_F=\bigcup \Gamma$, where the union is extended to every
cone $\Gamma$ such that $F$ is degenerate elliptic in $\Gamma$.
When we say that $F$ is degenerate elliptic, we mean that
$F$ is degenerate elliptic in $\Gamma_F\neq \emptyset$. A function $u\in C^2(\Omega)$ is said {\em admissible} for $F$ in $\Omega$ if $D^2u(x)\in\Gamma_F$ for every $x\in\Omega$.  In general, if not otherwise specified, we will consider for simplicity only operators such that $\Gamma_F=S_n$ throughout (then every regular function is admissible).

Given two functions $u$ and $\phi$ defined in an open set $\Omega$,
we say that $\phi$ \emph{touches $u$ by above at $x_0\in\Omega$} if
$$
\phi(x_0)=u(x_0)\textrm{ and }\phi(x)\geq u(x)\textrm{ in a neighbourhood of }x_0.
$$
Analogously,
we say that $\phi$ \emph{touches $u$ by below at $x_0\in\Omega$} if
$$
\phi(x_0)=u(x_0)\textrm{ and }\phi(x)\leq u(x)\textrm{ in a neighbourhood of }x_0.
$$
An upper semicontinuous function $u$ is a \emph{viscosity subsolution} of the equation
$F=0$ in $\Omega$ if, for every $C^2$ function $\phi$ touching $u$ by above at any point $x\in\Omega$, it holds
\begin{equation}\label{viscsubs}
F(x,u(x),D\phi(x),D^2\phi(x))\geq 0.
\end{equation}
A lower semicontinuous function $u$ is a \emph{viscosity supersolution} of $F=0$ in $\Omega$
if, for every admissible $C^2$ function $\phi$ touching $u$ by below at any point $x\in\Omega$, it holds
$$
F(x,u(x),D\phi(x),D^2\phi(x))\leq 0.
$$
A \emph{viscosity solution} is a continuous function which is
a viscosity subsolution and supersolution
of $F=0$ at the same time.

The technique proposed in this paper requires the use of the
comparison principle for viscosity solutions.
Since we will have only to compare a viscosity subsolution with a classical solution, we will need only a weak version of the comparison principle;
precisely,
we say that the operator $F$ satisfies \emph{the Comparison Principle} if the following statement holds:
\medskip

\noindent(CP)
{\em Let $u\in C(\overline{\Omega})\cap C^2(\Omega)$ and
$v\in C(\overline{\Omega})$ be respectively a classical supersolution
and a viscosity subsolution of $F=0$ such that $u\geq v$ on
$\partial \Omega$. 
Then $u\geq v $ in $\Omega$.}
\medskip


Comparison Principles for viscosity solutions are an actual and deep field of investigation and it is out of our aims to give here an updated picture of the state of the art, then we just refer to  \cite{caffacabre, UG, Koike}. However, when one of the involved function is regular, the situation is much easier and (CP) is for instance satisfied if $F$ is strictly proper, in other words if it is strictly monotone with respect to $u$.

\subsection{Minkowski addition and support functions of convex sets}

The Minkowski sum of two subsets $A_0$ and $A_1$ of $\RR^n$ is simply defined as follows
$$
A_0+A_1=\{ x+y\,:\, x\in A_0,\, y\in A_1\}\,.
$$
Let $\mu\in(0,1)$; the Minkowski convex combination of $A_0$ and $A_1$ (with coefficient $\mu$) is given by
$$
A_\mu=(1-\mu)A_0+\mu A_1=\{(1-\mu)x_0+\mu x_1\,:\, x_0\in\ A_0,\, x_1\in A_1\}\,.
$$
The famous {\em Brunn-Minkowski inequality} states 
\begin{equation}\label{BMineq}
|A_\mu|^{1/n}\geq(1-\mu)|A_0|^{1/n}+\mu|A_1|^{1/n}
\end{equation}
for every couple $A_0$, $A_1$ of measurable sets such that $A_\mu$ is also measurable. In other words, (\ref{BMineq}) states that the $n$-dimensional volume (i.e. Lebesgue measure) raised to power $1/n$ is concave with respect to Minkowski addition (see the beautiful paper by Gardner  \cite{Gardner} for a survey on this and related inequalities). 

When the involved sets are convex, Minkowski addition can be conveniently expressed in terms of support functions (see property (ii) below). 
The {\em support function} 
$h_\Omega:\RR^n\to\RR$ of a bounded convex set
$\Omega$ is defined as follows
$$
h_\Omega(X)=\max_{y\in\overline\Omega}\langle X,y\rangle\qquad X\in\RR^n\,.
$$
Every support function is convex and positively homogeneous of degree $1$, that is:
$$
h_\Omega(X+Y)\leq h_\Omega(X)+h_\Omega(Y)\,\,\mbox{ for every }X,Y\in\RR^n
$$
and
$$h_\Omega(tX)=t\,h_\Omega(X)\,\,\mbox{ for every }X\in\RR^n\mbox{ and }t\geq 0\,.
$$
Conversely, every convex and positively $1$-homogeneous function is the support function of a convex body  (i.e. a closed bounded convex set). This establishes
a one to one correspondence between support functions and convex bodies. 

Moreover the following properties hold:

(i) $h_{t\Omega}=th_\Omega$ for $t\geq0$;

(ii) $h_{\Omega_1+\Omega_2}=h_{\Omega_1}+h_{\Omega_2}$.
\\
The latter simply reads that the Minkowski addition of convex sets corresponds to the sum of support functions. 

As already said in the introduction, we denote the mean width of $\Omega$ by $w(\Omega)$, that is
$$w(\Omega)=\frac{1}{n\omega_{n}} \int_{S^{n-1}} \! \big(h(\Omega,\xi)+h(\Omega,-\xi)\big) \, d\xi=\frac{2}{n\omega_{n}} \int_{S^{n-1}} \! h(\Omega,\xi) \, d\xi\,.
$$
When $\Omega$ is a ball, $w(\Omega)$ coincides with its diameter. In the plane $w(\Omega)$ coincides with the perimeter of $\Omega$, up to a factor $\pi^{-1}$.


Given a convex set $\Omega$ and a point $x\in\partial\Omega$, we denote by $\nu_\Omega(x)$ the {\em exterior normal cone of $\Omega$ at $x$}, that is
$$
\nu_\Omega(x)=\{p\in\rn\,:\,\langle y-x,p\rangle\leq 0\,\,\mbox{for every }y\in\Omega\}\,.
$$
The normal cone of a convex set is a non-empty convex cone for every boundary point and in fact $\Omega$ is convex if and only if $\nu_\Omega(x)\neq\emptyset$ for every $x\in\partial\Omega$. The following elementary lemma about Minkowski addition will be useful in the sequel.
\begin{lemma}\label{Madditionlemma}
Let $\Omega_0,\,\Omega_1\subseteq\rn$ be open bounded convex sets and $\mu\in(0,1)$. 

Then $\Omega_\mu=(1-\mu)\Omega_0+\mu\Omega_1$ is an open bounded convex set; moreover
if $x_0\in\overline\Omega_0$ and $x_1\in\overline\Omega_1$ are such that $x=(1-\mu)x_0+\mu x_1\in\partial\Omega$,  then
$x_0\in\partial\Omega_0$, $x_1\in\partial\Omega_1$ and $\nu_{\Omega_\mu}(x)=\nu_{\Omega_0}(x_0)\cap\nu_{\Omega_1}(x_1)\neq\emptyset$.
\end{lemma}
The properties stated in the lemma can be considered folklore in the theory of convex bodies and the proof is straightforward.

For further details on convex sets, Minkowski addition and support functions, we refer to \cite{schneider}.

\subsection{Power concave functions}

Let $p\in[-\infty,+\infty]$ and $\mu\in(0,1)$. Given two real numbers $a>0$ and $b>0$, the quantity
\begin{equation}\label{pmedia}
 M_p(a,b;\mu)=\left\{\begin{array}{ll}
 \max\{a,b\} & p=+\infty
\\&\\
\left[(1-\mu)a^p+\mu b^p\right]^{1/p}
& \mbox{for }p\neq -\infty,\,0,\,+\infty\\ & \\
 a^{1-\mu}b^{\mu} & p=0\\&\\
\min\{a,b\}&p=-\infty
\end{array}\right.\end{equation}
is the {\em ($\mu$-weighted) $p$-mean} of $a$ and $b$.
For $a,b\geq 0$, we define $M_p(a,b;\mu)$ as above if $p\geq0$ and we set $M_p(a,b;\mu)=0$ if $p<0$ and $ab=0$.
Notice that $M_p$ is continuous with respect to $(a,b)\in[0,\infty)\times[0,\infty)$ for every $p$. See \cite{HLP} for more details.

A simple consequence of Jensen's inequality is that
\begin{equation}\label{disuguaglianzamedie}
M_p(a,b;\mu)\leq M_q(a,b;\mu)\quad \mbox{if }-\infty\leq p\leq
q\leq+\infty\,.\end{equation}

\begin{definition}
\label{Definition:1.1}
Let $\Omega$ be an open convex set in $\RR^n$ and $p\in[-\infty,\infty]$.
A function $v:\Omega\to[0,+\infty)$ is said {\em $p$\,-concave} if
$$
v((1-\mu)x+\mu y)\geq M_p(v(x),v(y);\mu)
$$
for all $x$, $y\in \Omega$ and $\mu\in(0,1)$. 

In the cases $p=0$ and $p=-\infty$, 
$v$ is also said {\em log-concave} and {\em quasi-concave} in $\Omega$, respectively. 
\end{definition}
In other words, a non-negative function $v$, with convex support $\Omega$, is $p$-concave if:\\
- it is a non-negative constant in $\Omega$, for $p=+\infty$;\\
- $v^p$ is concave in $\Omega$, for $p>0$;\\
- $\log v$ is concave in $\Omega$, for $p=0$;\\
- $v^p$ is convex in $\Omega$, for $p<0$;\\
- it is quasi-concave, i.e. all of its superlevel sets are convex, for $p=-\infty$.
\\
Notice that $p=1$ corresponds to usual concavity

It follows from (\ref{disuguaglianzamedie}) that 
{\em if $v$ is $p$\,-concave, 
then $v$ is $q$\,-concave for any $q\le p$}. 
Hence quasi-concavity is the weakest conceivable concavity property. 
\vspace{5pt} 

It is well known that solutions of elliptic Dirichlet problems in convex domains are often power concave. For instance, a famous result by Brascamp and Lieb \cite{bl} says that the first positive eigenfunction of the Laplace operator in a convex domain is log-concave; another classical result states that the square root of the solution to the torsion problem in a convex domain is concave, see \cite{Kawohl, Kenni, ML}. Power concave solutions have been also studied in
\cite{kaw,Kaw00, Kore2} and more recent developments are for instance in \cite{all, Juutinen, LeeVazquez, LiuMaXu, MaXu, Salani, Ye}; furthermore see \cite{cuoghis} and  \cite{BS}, which are strongly related to the present paper.

\subsection{The Borell-Brascamp-Lieb inequality}

The Borell-Brascamp-Lieb inequality (see \cite{borel, bl})  is a generalization of the Pr\'ekopa-Leindler inequality. 
I recall it here in the form taken from \cite[Theorem 10.1]{Gardner}. 
\begin{proposition}
\label{BBLprop}
Let $\mu\in(0,1)$, $f,g,h$ nonnegative functions in $L^1(\RR^n)$, 
and $-1/n\leq s\leq \infty$. 
Assume that
\begin{equation}
\label{eq:4.2}
h\big((1-\mu)x+\mu y\big)\geq M_s(f(x),g(y);\mu)
\end{equation}
for all $x\in\mbox{\em sprt}(f),\,y\in\mbox{\em sprt}(g)$. Then
$$
\int_{\RR^n} h\,dx\ge M_q\left(\int_{\RR^n} f\,dx,\int_{\RR^n} g\,dx\,;\mu\right)\,,
$$
where 
\begin{equation}
\label{eq:4.1}
q=\left\{
\begin{array}{ll}
1/n &\mbox{if}\quad s=+\infty,\vspace{5pt}\\
s/(ns+1) &\mbox{if}\quad s\in(-1/n,+\infty),\vspace{5pt}\\
-\infty&\mbox{if}\quad s=-1/n.
\end{array}
\right.
\end{equation}
\end{proposition}

The Pr\'ekopa-Leindler inequality corresponds to the case $s=0$ and it is a functional version of the Brunn-Minkowski inequality.

\section{The $(p,\mu)$-convolution of non-negative functions}

\vskip 0,3 cm From now on, throughout the paper,  we consider two open bounded convex sets $\Omega_0,\,\Omega_1\subset\rn$ and a fixed real number $\mu\in(0,1)$, and denote by $\Omega_\mu$ the Minkowski convex combination (with coefficient $\mu$) of $\Omega_0$ and $\Omega_1$,
i.e. $\Omega_\mu=(1-\mu)\Omega_0+\mu\Omega_1$.
\begin{definition}\label{defupmu}
Let $p\in\RR$, $\mu\in(0,1)$, $u_0\in C(\overline\Omega_0)$ and $u_1\in C(\overline\Omega_1)$
such that $u_i\geq 0$ in $\overline\Omega_i$, $i=0,1$.
The {\em $(p,\mu)$-convolution}
of $u_0$ and $u_1$ is the function $u_{p,\mu}:\overline \Omega_\mu \rightarrow
\RR$ defined as follows:
\begin{equation}\label{pconcaveenvelope}
\begin{array}{rl}u_{p,\mu}(x)=
\sup\big\{&\!\!\!\!\!M_p\big(u_0(x_0),u_1(x_1);\mu\big)\,:\,\\
&\,\,\,x=(1-\mu)x_0+\mu x_1\,,\,x_i\in\overline{\Omega_i},\,i=0,1\big\}.
\end{array}
\end{equation}
\end{definition}
The above definition can be extended to the case $p=\pm\infty$, but we do not need here. Let me recall however that the case $p=-\infty$ has been useful in \cite{Bo84, CoSa} to prove the Brunn-Minkowski inequality for $p$-capacity of convex sets.

Let $p\neq 0$; then, roughly speaking, the graph of $u_{p,\mu}^p$ is obtained as the Minkowski convex combination (with coefficient $\mu$) of the graphs of $u_0^p$ and $u_1^p$; precisely we have
$$
K^{(p)}_{\mu}=(1-\mu)K^{(p)}_0+\mu K^{(p)}_1\,,
$$
where
$$
K^{(p)}_{\mu}=\{(x,t)\in\RR^{n+1}\,:\, x\in\Omega_\mu,\,0\leq t\leq u_{p,\mu}(x)^p\}\,,
$$
$$
K^{(p)}_i=\{(x,t)\in\RR^{n+1}\,:\, x\in\Omega_i,\,0\leq t\leq u_i(x)^p\}\,,\quad i=0,1\,.
$$
In other words, the $(p,\mu)$-convolution of $u_0$ and $u_1$ corresponds to the $(1/p)$-power of the supremal convolution (with coefficient $\mu$) of $u_0^p$ and $u_1^p$. When $p=0$, the above geometric considerations continue to hold with logarithm in place of power $p$ and exponential in place of power $1/p$. When $p=1$, $u_{1,\mu}$ is just the usual supremal convolution of $u_0$ and $u_1$.
For more details on infimal/supremal convolutions of convex/concave functions, see \cite{rock, st} (and also  \cite{CoCuSa, salanima}).

From Definition \ref{defupmu} and (\ref{disuguaglianzamedie}), for every $\mu\in(0,1)$ we get
\begin{equation}\label{disuguaglianzapenvelope}u\leq u_{p,\mu}\leq u_{q,\mu}\quad \mbox{for } -\infty\leq p\leq q\leq+\infty\,.\end{equation}
Clearly $u$ is $p$-concave if and only if $u=u_{p,\mu}$ for every $\mu\in(0,1)$.



\begin{lemma} \label{proprietaup} Let $p\in[-\infty,+\infty)$, $\mu\in(0,1)$. For $i=0,1$ let $u_i\in C(\overline \Omega_i)$ such
that $u_i=0$ on $\partial \Omega_i$ and $u_i>0$ in $\Omega_i$.
Then $u_{p,\mu}\in C(\overline
\Omega_\mu)$ and
\begin{equation} \label{bordoup}u_{p,\mu}>0
\quad \mbox{in } \Omega_\mu,\quad u_{p,\mu}=0\quad \mbox{on }
\partial \Omega_\mu.\end{equation}

\end{lemma}
{\em Proof. } The proof of this lemma is almost straightforward and completely analogous to the proof of
 \cite[Lemma 1]{BS}. We just notice that $u_{p,\mu}>0$ in $\Omega$  by the very definition of $u_{p,\mu}$ while $u_{p,\mu}=0$  on $\partial\Omega$ by Lemma \ref{Madditionlemma}.
\qed\medskip

Notice that, as $\overline \Omega_i$ is compact for $i=0,1$ and $M_p$, $u_0$ and $u_1$ are
continuous, then the supremum in (\ref{pconcaveenvelope}) is in fact a maximum.
Hence for every $\bar{x}\in\overline\Omega_\mu$ there exist
$x_0\in \overline \Omega_0$ and $x_1\in\overline\Omega_1$
such that
\begin{equation} \label{massimo0}
\bar{x}=(1-\mu) x_0+\mu x_1\,,\qquad
u_{p,\mu}(\bar{x})=M_p(u_0(x_0),u_1(x_1);\mu)\,.
\end{equation}

The next lemma is fundamental to this paper.

\begin{lemma}
\label{nobordo_p_maggiore0}
Let $p\in[0,1)$, $\mu \in(0,1)$, 
$u_i\in C^1(\Omega_i)\cap C(\overline \Omega_i)$ such
that $u_i=0$ on $\partial \Omega_i$, $u_i>0$ in $\Omega_i$ for $i=0,1$.

In case $p>0$ assume furthermore that for $i=0,1$ it holds
\begin{equation}
\label{liminf}
\liminf_{y\rightarrow x}
\frac{\partial u_i(y)}{\partial \nu}>0
\end{equation}
for every $x\in\partial\Omega_i$, where $\nu$ is any inward direction of $\Omega_i$ at $x$.

If $\bar{x}$ lies in the interior of $\Omega_\mu$,
then the
points $x_{0}$ and $x_{1}$ defined by (\ref{massimo0}) belong to the interior of $\Omega_0$ and $\Omega_1$, respectively, and
\begin{equation} \label{moltiplicatorilagrange} u_0(x_{0})^{p-1}D
u_0(x_{0})=u_1(x_{1})^{p-1}Du_1(x_{1})\,.
\end{equation}
\end{lemma}
{\em Proof. }
First we prove that $x_i\in\Omega_i$ for $i=0,1$.

The case $p=0$ easily follows from (\ref{bordoup}) and the definition of $M_0$, since $u_{p,\mu}(\bar{x})>0$ while
$u_0(x_0)^{1-\mu}u_1(x_1)^\mu=0$ if $x_0\in\partial\Omega_0$ or $x_1\in\partial\Omega_1$.

Then let $p>0$. By contradiction, assume that (up to a relabeling) $x_0\in\partial\Omega_0$. Then $u_0(x_0)=0$ and $x_1$ must lie in the interior of 
$\Omega_1$, otherwise
 $u_{p,\mu}(\bar{x})=0$, contradicting (\ref{bordoup}). Notice that in this case
$$
u_{p,\mu}(\bar{x})=\mu^{1/p}u_1(x_1)\,.
$$
Set $v_0=u_0^p$, $v_1=u_1^p$ and
$$a=|Dv_1(x_{1})|=p\,u_1(x_{1})^{p-1}|Du_1(x_{1})|\,.
$$
By the regularity of $u_1$, we have
\begin{equation}\label{neweq1}
|Dv_1|<a+1\quad\mbox{in }B(x_{1},r_{1})\subset\Omega_1
\end{equation}
for $r_{1}>0$ small enough.

Now take any direction $\nu$
pointing inwards into $\Omega_0$ at $x_0$;
by assumption (\ref{liminf}) we get 
\begin{equation}\label{gradinfty}
\liminf_{x\rightarrow x_0}\,
\frac{\partial v_0(x)}{\partial \nu}=+\infty\,,
\end{equation}
whence
\begin{equation}\label{neweq2}
\frac{\partial v_0}{\partial\nu}>a+1\quad\mbox{in }\Omega_0\cap B(x_0,r_0)
\end{equation}
for $r_0>0$ small enough.

Next we take $\rho<\min\{(1-\mu) r_0,\,\mu r_{1}\}$ and
we consider the points
$$
\begin{array}{ll}
\tilde x_0=x_0+\frac{\rho}{(1-\mu)} \nu\,,\\
\\
\tilde x_{1}=x_{1}-\frac{\rho}{\mu}\nu\,.
\end{array}
$$
We have
$$
\tilde x_0\in B(x_0,r_0)\cap\Omega_0\,,\quad \tilde x_{1}\in B(x_{1},r_{1})
$$
and
\begin{equation}\label{xxtilda}
\bar{x}=(1-\mu)\tilde x_0+\mu\tilde x_1\,.
\end{equation}
Then from (\ref{neweq1}) and (\ref{neweq2}) we get
\begin{eqnarray*}
&u_0(\tilde x_0)^p=v_0(\tilde x_0)>v_0(x_0)+(a+1)\frac{\rho}{(1-\mu)}=(a+1)\frac{\rho}{(1-\mu)}\,,\\
&u_1(\tilde x_{1})^p=v_1(\tilde x_{1})\geq v_1(x_{1})^p-(a+1)\frac{\rho}{\mu}=u_1(x_{1})^p-(a+1)\frac{\rho}{\mu}\,,
\end{eqnarray*}
whence
$$
\begin{array}{rl}
\big[(1-\mu) u_0(\tilde x_0)^p&\!\!\!\!+\,\mu\,u_1(\tilde x_1)^p\big]^{1/p}>\\
&\,\,\left[(1-\mu)(a+1)\frac{\rho}{(1-\mu)}
+\mu u_1(x_{1})^p-\mu  (a+1)\frac{\rho}{\mu}\right]^{1/p}=u_{p,\mu}(\bar{x})
\end{array}
$$
which contradicts the definition of $u_{p,\mu}$, due to (\ref{xxtilda}).

So far, we have proved that $x_i$ must stay in the interior of $\Omega_i$ for $i=0,1$. Then by the Lagrange
Multipliers Theorem we easily get (\ref{moltiplicatorilagrange}). In fact, we could just notice that $x_0\in\Omega_0$ is an interior maximum point for the function
$$
f(x)=M_p\left(u_0(x),u_1\big(\frac{\bar{x}-(1-\mu)x}{\mu}\big);\mu\right)
$$
and $\nabla f(x_0)=0$ gives (\ref{moltiplicatorilagrange}).

The proof of the lemma is complete.
\qed\medskip

\subsection{The $(p,\mu)$-convolution of more than two functions.}
The definition of the $(p,\mu)$-convolution of two functions is easily extended to an arbitrary number of functions.  

Let $3\leq m\in\NN$ and set $\Gamma_m^+=\{(x_1,\dots,x_m)\in\RR^m\,:\,x_i\geq 0\,,\,i=1,\dots,m\}$ and
$$\Gamma_m^1=\left\{(\mu_1,\dots,\mu_m)\in\Gamma_m^+\,:\,\mu_i>0\mbox{ for }i=1,\dots,m\mbox{ and }\sum_{i=1}^m\mu_i=1\right\}\,.
$$
Let $p\in[-\infty,+\infty]$, $\mu\in\Gamma_m^1$ and $a=(a_1,\dots,a_m)\in\Gamma_m^+$. 
If $\prod_{i=1}^ma_i>0$, the $p$-mean of $a_1,\dots,a_m$ with coefficient $\mu$ is defined as follows:
$$
 M_p(a_1,\dots,a_m;\mu)=\left\{\begin{array}{ll}
 \max\{a_1,\dots,a_m\} & p=+\infty
\\&\\
\left[\sum_{i=1}^m\mu_ia_i^p\right]^{1/p}
& p\neq -\infty,\,0,\,+\infty\\ & \\
\prod_{i=1}^ma_i^{\mu_i}& p=0\\&\\
\min\{a_1,\dots,a_m\}&p=-\infty
\end{array}\right.$$
If $\prod_{i=1}^ma_i=0$, we define $M_p(a,\mu)$ as above if $p\geq0$ and we set $M_p(a,\mu)=0$ if $p<0$.

If we now consider $m$ non-negative
functions $u_1,u_2,\dots,u_m$ supported in the sets $\Omega_1,\Omega_2,\dots,\Omega_m$ respectively,  we can define
\begin{equation}\label{pconcaveenvelopem}
\begin{array}{rl}u_{p,\mu}(x)=
\sup\big\{&\!\!\!\!M_p
\left(u_0(x_0),\dots,u_m(x_m);\mu\right)\,:\\
&\,x_i\in\overline{\Omega_i},\,i=1,\dots,m,\,
\,\,\,\,\,x=\sum_{i=1}^m\mu_ix_i\big\}.
\end{array}
\end{equation}
Clearly all the properties and lemmas stated and proved before for the case $m=2$ continue to hold in the case $m\geq 3$, with the obvious modifications.
In particular we explicitly write the following.

\begin{lemma} \label{proprietaupm} Let $p\in[-\infty,+\infty)$,  $\mu\in\Gamma_m^1$. Let $u_i\in C(\overline \Omega_i)$ such
that $u_i=0$ on $\partial \Omega_i$ and $u_i>0$ in $\Omega_i$, for $i=1,\dots,m$.
Then $u_{p,\mu}\in C(\overline
\Omega_\mu)$ and
\begin{equation} \label{bordoupm}u_{p,\mu}>0
\quad \mbox{in } \Omega_\mu,\quad u_{p,\mu}=0\quad \mbox{on }
\partial \Omega_\mu.\end{equation}
\end{lemma}

As before, since $\overline \Omega_i$ is compact for $i=1,\dots, m$ and $M_p$, $u_1,\dots,u_m$ are
continuous, the supremum in (\ref{pconcaveenvelopem}) is in fact a maximum.
Hence for every $\bar{x}\in\overline\Omega_\mu$ there exist
$x_0\in \overline \Omega_0,\dots,x_m\in\overline\Omega_m$
such that
\begin{equation} \label{massimo0m}
\bar{x}=\sum_{i=1}^m\mu_ix_i\,,\qquad
u_{p,\mu}(\bar{x})=M_p(u_1(x_1),\dots,u_m(x_m);\mu)\,.
\end{equation}
\begin{lemma}
\label{nobordo_p_maggiore0m}
Let $p\in[0,1)$, $\mu \in(0,1)$, 
$u_i\in C^1(\Omega_i)\cap C(\overline \Omega_i)$ such
that $u_i=0$ on $\partial \Omega_i$, $u_i>0$ in $\Omega_i$ for $i=1,\dots,m$.
In case $p>0$ assume furthermore that for  (\ref{liminf}) holds for $i=1,\dots,m$.

If $\bar{x}$ lies in the interior of $\Omega_\mu$,
then the
points $x_1,\dots,x_m$ defined by (\ref{massimo0}) belong to the interior of $\Omega_1,\dots,\Omega_m$, respectively, and
\begin{equation} \label{moltiplicatorilagrangem} u_1(x_{1})^{p-1}D
u_1(x_{1})=\dots=u_m(x_{m})^{p-1}Du_m(x_{m})\,.
\end{equation}
\end{lemma}

\section{The main theorem}
As before and throughout, $\Omega_0$ and $\Omega_1$ are open bounded convex sets in $\rn$, $\mu\in(0,1)$ and $\Omega_\mu=(1-\mu)\Omega_0+\mu\Omega_1$.

For $i=0,1,\mu$, we denote by $u_i$ a solution of the following problem
$$
(P_i)\quad
\left\{
\begin{array}{ll}
F_i(x,u_i,D u_i, D^2u_i)=0 & \textrm{ in } \Omega_i\,,\\
u_i=0 & \textrm{ on } \partial \Omega_i\,,\\
u_i>0 &\textrm{ in }\Omega_i\,,
\end{array}
\right.
$$
where $F_i:\Omega_i\times
[0,+\infty)\times\RR^n\times S_n$ is a proper elliptic operator.

If not otherwise specified, we will consider {\em classical solutions for $i=0,1$} (that is: $u_0\in C^2(\Omega_0)\cap C(\overline\Omega_0)$ and
$u_1\in C^2(\Omega_1)\cap C(\overline\Omega_1)$ 
and they satisfy pointwise everywhere all the equations in $(P_0)$ and $(P_1)$, respectively), while $u_\mu\in C(\overline\Omega_\mu)$ may be a {\em viscosity solution} of the corresponding problem $(P_\mu)$.

For $i=0,1,\mu$ and for every fixed $(\theta,p)\in\RR^n\times[0,\infty)$  we define
$G_{i,p}^{(\theta)}:
\Omega_i\times (0,+\infty)\times S_n \rightarrow \RR$ as
\begin{equation}\label{Gpno0}
G_{i,p}^{(\theta)}(x,t,A)= F_i(x,t^{\frac{1}{p}},t^{\frac{1}{p}-1}\theta,
t^{\frac{1}{p}-3}A)\,\quad\mbox{ for }p>0\,,
\end{equation}
and
\begin{equation}\label{G0}
G_{i,0}^{(\theta)}(x,t,A)=F_i(x,e^t,e^t\theta,
e^tA)\,.
\end{equation}

{\bf{Assumption $(A_{\mu,p})$.}}
{\em Let $\mu\in(0,1)$ and $p\geq0$. We say that} $F_0,F_1,F_\mu$ satisfy the assumption $(A_{\mu,p})$ {\em if, for every fixed $\theta\in\rn$, the following holds:
$$
\begin{array}{ll}
G^{(\theta)}_{\mu,p}\big((1-\mu)x_0+\mu x_1, &\!\!\!\!(1-\mu)t_0+\mu t_1,(1-\mu)A_0+\mu A_1\big)\geq\\
&\qquad\min\{
G_{0,p}^{(\theta)}(x_0,t_0,A_0);\,G_{1,p}^{(\theta)}(x_1,t_1,A_1)\}
\end{array}
$$
for every $x_0\in\Omega_0$, $x_1\in\Omega_1$, $t_0,t_1>0$ and $A_0,A_1\in S_n$.}
\medskip

Now we are ready to state the main result of the paper.

\begin{theorem}\label{mainthm}
Let $\mu\in(0,1)$ and $\Omega_i$ and $u_i$, $i=0,1,\mu$, be as above described.
Assume that the operator $F_\mu$ satisfies the comparison principle $(CP)$ and that $F_0,F_1,F_\mu$ satisfy the assumption $(A_{\mu,p})$ for some $p\in[0,1)$.
If $p>0$, assume furthermore that (\ref{liminf}) holds true for $i=0,1$.

Then
\begin{equation}\label{eqmainthm}
u_\mu((1-\mu)x_0+\mu\,x_1)\geq M_p(u_0(x_0),u_1(x_1);\mu)
\end{equation}
for every $x_0\in\Omega_0,\,x_1\in\Omega_1$.
\end{theorem}

We remark that assumption (\ref{liminf}) is not needed for $p=0$, while for $p>0$ it is in general provided by a suitable version of the Hopf's Lemma. Notice also that,  for $p<1$, (\ref{liminf}) implies (\ref{gradinfty}).
In fact, we could also apply our argument to the case $p\geq1$; in such a case however we would need to assume directly (\ref{gradinfty}) instead of 
(\ref{liminf}).

\medskip
Coupling (\ref{eqmainthm}) with the Borell-Brascamp-Lieb inequality (i.e. Proposition \ref{BBLprop}) leads to a comparison of the $L^r$ norms of $u_\mu$ with suitable combinations of the $L^r$ norms of $u_0$ and $u_1$. Precisely, we have the following corollary.
\begin{corollary}\label{maincor}
In the same assumptions and notation of Theorem \ref{mainthm}, for every $r>0$ we have
\begin{equation}\label{eqLpnorm}
\|u_\mu\|_{L^r(\Omega_\mu)}\geq M_q(\|u_0\|_{L^r(\Omega_0)},\|u_1\|_{L^r(\Omega_1)};\mu)\,,
\end{equation}
where 
$$
q=\left\{\begin{array}{ll}\frac{pr}{np+r}\qquad&\mbox{for }r\in(0,+\infty)\\
\\
p\qquad&\mbox{for }r=+\infty\,.
\end{array}\right.
$$
\end{corollary}
{\em Proof. }
The inequality for the $L^\infty$ norms is a straightforward consequence of  (\ref{eqmainthm}), obtained by taking $x_0$ and $x_1$ as points which realize the maximum of $u_0$ and $u_1$, respectively (in fact in this case equality holds in (\ref{eqLpnorm})). The proof of the inequality for a generic $r\in(0+\infty)$ follows from Proposition \ref{BBLprop}, applied to the functions $h=u_\mu^r$, $f=u_0^r$ and $g=u_1^r$ with $s=p/r$, assumption (\ref{eq:4.2}) being satisfied thanks to (\ref{eqmainthm}).
\qed\medskip

Notice that in some special cases, involving particular operators, results similar to those we could obtain by applying Theorem \ref{mainthm} and Corollary \ref{maincor} to the situations at hands, has been already proved (even though not explicitly stated) and used to prove Brunn-Minkowski type inequalities for variational functionals, see for instance \cite{Co, CoCuSa, LiuMaXu, Salani, WX}. Indeed, Theorem \ref{mainthm} could be regarded as a general Brunn-Minkowski inequality for solutions of PDE's (and then applied to obtain Brunn-Minkowskii type inequalities for possibly related functionals).


\section{Proof of Theorem \ref{mainthm}}

The proof of Theorem 1.1 essentially consists of the following lemma.

\begin{lemma}\label{lemmapleq0}
In the same assumptions and notation of Theorem \ref{mainthm}, it follows that $u_{p,\mu}$ is a viscosity subsolution of problem $(P_\mu)$.
\end{lemma}
{\em Proof. }
The proof follows somehow the steps of \cite{BS, cuoghis, IS5} and the strategy is the following: for every $\bar{x}
\in\Omega_\mu$, we construct a function $\varphi_{p,\mu}\in C^2(\Omega_\mu)$
which touches $u_{p,\mu}$
by below at $\bar{x}$ and such that
\begin{equation}
F(\bar{x},\varphi_{p,\mu}(\bar{x}),
D \varphi_{p,\mu}(\bar{x}),
D^2 \varphi_{p,\mu}(\bar{x}))\geq 0.
\label{Fphileq0}\end{equation}
Clearly this implies that $u_{p,\mu}$
is a viscosity subsolution of $(P_\mu)$: indeed every test function $\phi$ touching $u_{p,\mu}$ at $\bar{x}$ by above must also touch $\varphi_{p,\mu}$ at $\bar{x}$ by above, then $$\phi(\bar{x})=\varphi_{p,\mu}(\bar{x})\,,\,\,D\phi(\bar{x})=D\varphi_{p,\mu}(\bar{x})\,\,\mbox{ and }\,\,D^2\phi(\bar{x})\geq D^2\varphi_{p,\mu}(\bar{x})$$
and (\ref{viscsubs}) follows from the ellipticity of $F$.
\medskip

Then consider $\bar{x} \in \Omega$. By
Lemma \ref{nobordo_p_maggiore0}, there exist $x_{0}\in\Omega_0$ and $x_{1}\in
\Omega_1$ satisfying (\ref{massimo0}) and such that (\ref{moltiplicatorilagrange}) holds.

First we treat the case $p>0$ and, for a  small enough
$r>0$, we introduce the function
$\varphi_{p,\mu}:B(\x, r) \rightarrow \RR$ defined as follows
\begin{equation}\label{approssimantesotto}\varphi_{p,\mu}(x)=
\big[(1-\mu)u_0\left(x_{0}+a_{0}(x-\x)\right)^p+\mu u_1\left(x_{1}+a_{1}(x-\x)\right)^p
\big]^{1/p}\end{equation} where
\begin{equation}\label{defai}a_{i}=\frac{u_i(x_{i})^p}{u_{p,\mu}(\x)^p},\quad\mbox{for
}i=0,1.\end{equation}

The following facts trivially hold:\begin{itemize}
    \item[(A)] $(1-\mu)a_{0}+\mu a_{1}=1$ by (\ref{massimo0});
    \item [(B)]$x=(1-\mu)(x_{0}+a_{0}(x-\x))+\mu(x_{1}+a_{1}(x-\x))$ for every
$x\in B(\x,r)$, thanks to (A) and the first equation in (\ref{massimo0});
    \item[(C)] $\varphi_{p,\mu}(\x)=u_{p,\mu}(\x)$;
    \item [(D)]$\varphi_{p,\mu}(x)\leq u_{p,\mu}(x)$ in $B(\x,r)$ (this
follows from $(B)$ and the definition of
    $u_{p,\mu}$).
\end{itemize}
In particular, (C) and (D) say that $\varphi_{p,\mu}$ touches $u_{p,\mu}$
from below at $\x$.

\vskip 0,5 cm A straightforward calculation yields
$$D \varphi_{p,\mu}
(\x)=\varphi_{p,\mu}(\x)^{1-p}\left[(1-\mu)u_0(x_{0})^{p-1}a_{0}\,D u_0(x_{0})+\mu u_1(x_{1})^{p-1}a_{1}\,D u_1(x_{1})\right],
$$
Then,
by (\ref{moltiplicatorilagrange}), (\ref{defai}) and the definition of
$\varphi_{p,\mu}$, we get
\begin{equation}D \varphi_{p,\mu}(\x)=\varphi_{p,\mu}(\x)^{1-p}u_i(x_{i,p})^{p-1}D u_i(x_{i,p})
\textrm{ for }i=0,1.\label{gradientivarphi}
\end{equation}
Thanks to another straightforward calculation and using (\ref{moltiplicatorilagrange}), (\ref{defai}), (\ref{gradientivarphi}) and the definition of
$\varphi_{p,\mu}$, we also obtain
\begin{eqnarray*} D^2\varphi_{p,\mu}(\x)&=& (1-\mu)\frac{u_0(x_{0})^{3p-1}}{\varphi_{p,\mu}(\x)^{3p-1}}D^2u_0(x_{0})+
\mu\frac{u_1(x_{1})^{3p-1}}{\varphi_{p,\mu}(\x)^{3p-1}}D^2u_1(x_{1})+\\&&\\&&
+(1-p)\varphi_{p,\mu}(\x)^{-1}
A\,D\varphi_{p,\mu}(\x)\otimes D \varphi_{p,\mu}(\x)\,,\end{eqnarray*}
where
$$
A=1-\varphi_{p,\mu}(\x)^{-p}[(1-\mu)u_0(x_{0})^p+\mu u_1(x_{1})^p]\,.
$$
Now notice that (C) and (\ref{massimo0}) give
$$
A=0\,.
$$
Then
\begin{equation} D^2\varphi_{p,\mu}(\x)=(1-\mu)\frac{u_0(x_{0})^{3p-1}}{\varphi_{p,\mu}(\x)^{3p-1}}D^2u_0(x_{0})+
\mu\frac{u_1(x_{1})^{3p-1}}{\varphi_{p,\mu}(\x)^{3p-1}}D^2u_1(x_{1})\,
.\label{hessianafinale}\end{equation}
\\
Since $u_0$ and $u_1$ are classical solutions of $(P_0)$ and $(P_1)$ respectively, it
follows that for $i=0,1$
$$
G^{(\theta)}_{i,p}(x_{i},u_i(x_{i})^p,
u_i(x_{i})^{3p-1}D^2u_i(x_{i}))=F_i(x_{i},u_i(x_{i}),Du_i(x_{i}),
D^2u_i(x_{i}))=0\,,
$$ 
where
$$
\theta=\varphi_{p,\mu}(\x)^{p-1}
D \varphi_{p,\mu}(\x)\,.
$$
Then, by setting  $\mu_0=(1-\mu)$ and $\mu_1=\mu$, assumption $(A_{\mu,p})$ entails
$$
G^{(\theta)}_{\mu,p}\left(\sum_{i=0}^{1}\mu_i
x_{i},
\sum_{i=0}^{1}\mu_iu_i(x_{i})^p,
\sum_{i=0}^{1}
\mu_i\,u_i(x_{i})^{3p-1}D^2u_i(x_{i}) \right)\geq 0\,,$$
and thanks to $(C)$ and (\ref{hessianafinale}) this  precisely coincides with
$$
G^{(\theta)}_{\mu,p}(\x,\varphi_{p,\mu}(\x)^p,\varphi_{p,\mu}(\x)^{3p-1}D^2\varphi_{p,\mu}(\x))\geq 0\,.
$$
The latter implies (\ref{Fphileq0}) by the definition of $G^{(\theta)}_{\mu,p}$ and this concludes the proof for $p>0$.
\medskip

The case $p=0$ is similar, the only difference
consisting in that we set
$$\varphi_{0,\mu}:=\exp\big((1-\lambda) \log u_0(x_{1,0}+x-\x)+\mu\log u_1(x_{n+1,0}+x-\x  )\big)\,,$$ which
means $a_{i,0}=1$ for $i=0,1$.
\qed
\medskip

The proof of  Theorem \ref{mainthm} is now very easy.
\medskip

\noindent{\em Proof of  Theorem \ref{mainthm}.}
Under the assumptions of the theorem, we can apply the previous lemma
to obtain that $u_{p,\mu}$ is a viscosity subsolution of $(P_\mu)$.
Then by the Comparison Principle we get the thesis.
\qed\medskip

\subsection{A generalization.}
Looking at the proof of Lemma \ref{lemmapleq0}, it is easily understood that assumption $(A_{\mu,p})$ can be in fact substituted by a slightly weaker one: precisely what really matters is that the inequality in $(A_{\mu,p})$ holds only for $(x_i,t_i,A_i)$ such that $G_{p,\theta}(x_i,t_i,A_i)=0$, $i=0,1$.
Moreover, it is clear that, when considering the combination of more than two Dirichlet problems, a generalized version of Theorem \ref{mainthm} continues  to hold. 

Exactly: let $m\in\NN$, $m\geq 2$, and $\mu=(\mu_1,\mu_2,\dots,\mu_m)\in\Gamma^1_m$; let $\Omega_i$, $F_i$ and $u_i$ be a convex set, a proper elliptic operator and the solution of problem $(P_i)$, respectively, for $i=1,\dots,m$ and $i=\mu$, where
$$
\Omega_\mu=\sum_{i=1}^m\mu_i\Omega_i\,;
$$
define $G_{i,p}^{(\theta)}$ as in (\ref{Gpno0}) and (\ref{G0}) and set
$$
Z_{i,p}^{(\theta)}=\{(x,t,A)\,:\,G_{i,p}^{(\theta)}(x,t,A)=0\}
$$
for $i=1,\dots,m$;
then we say that the operators $F_\mu, F_1,\dots, F_m$ satisfies the  {\em Assumption Weak $(A_{\mu,p})$} if
$$
(WA_{\mu,p})\qquad S_{\mu,p}^{(\theta)}=\{(x,t,A)\,:\,G_{\mu,p}^{(\theta)}(x,t,A)\geq0\}\supseteq\sum_{i=1}^m\mu_iZ_{i,p}^\theta
$$
for every $\theta\in\RR^n$.

\begin{theorem}\label{mainthm2}
Assume that the operator $F_\mu$ satisfies the comparison principle $(CP)$ and that $F_1,\dots,F_m,F_\mu$ satisfy the assumption $(WA_{\mu,p})$ for some $p\in[0,1)$.

If $p>0$, assume furthermore that for $i=1\dots,m$ it holds
\begin{equation}
\label{liminf2}
\liminf_{y\rightarrow x}
\frac{\partial u_i(y)}{\partial \nu}>0
\end{equation}
for every $x\in\partial\Omega_i$, where $\nu$ is any inward direction of $\Omega_i$ at $x$.

Then
\begin{equation}\label{eqmainthm2}
u_\mu\big(\sum\mu_ix_i\big)\geq M_p(u_1(x_1),\dots,u_m(x_m);\mu)
\end{equation}
for every $x_1\in\Omega_1,\,x_2\in\Omega_2\dots,\,x_m\in\Omega_m$.
\end{theorem}

Obviously the key point is that it holds an appropriate version of Lemma \ref{lemmapleq0}, that is the following.
\begin{lemma}\label{lemmapleq02}
In the same assumptions and notation of Theorem \ref{mainthm2}, it follows that $u_{p,\mu}$, defined by (\ref{pconcaveenvelopem}), is a viscosity subsolution of (\ref{pb1}) in $\Omega_\mu$.
\end{lemma}

The proof of this lemma is just a straightforward adaptation of the proof of Lemma \ref{lemmapleq0} and we omit it.

\section{Rearrangements}\label{rearrsect}

Throughout $p$ will be a real positive number and $\Omega\subset\rn$ an open bounded convex set. 

We say that $\Omega^\sharp_m$ is a {\em rotation mean} of $\Omega$ if there exist a number $m\in\mathbb{N}$ and $\rho_1,\dots,\rho_m\in SO(n)$ such that
$$
\Omega^\sharp_m=\frac{1}{m}\left(\rho_1\Omega+\dots+\rho_m\Omega\right)\,.
$$
The following theorem is due to Hadwiger. 
\begin{theorem}{\cite[Theorem 3.3.2]{schneider}}.\label{Hadwiger} Given an open bounded convex set $\Omega$, there exists a sequence of rotation means of $\Omega$ converging in Hausdorff metric to a ball $\Omega^\sharp$ 
with 
diameter equal to the mean width $w(\Omega)$ of $\Omega$.
\end{theorem}

Let $u\in C(\overline\Omega)$ be a non-negative function, positive in $\Omega$ and vanishing on $\partial\Omega$,
and let $\rho\in SO(n)$; we set
\begin{equation}\label{urho}
u_\rho(x)=u(\rho^{-1}x)\quad\mbox{for } x\in\rho\overline\Omega\,.
\end{equation}
Now let $\{\rho_i\}_{i=1}^\infty\subset SO(n)$ be the sequence of rotations associated to $\Omega$ by Theorem \ref{Hadwiger}, such that $\Omega^\sharp_m$ converges to $\Omega^\sharp$, and set
$$\Omega_i=\rho_i\Omega\,,\quad u_i=u_{\rho_i}\,\,\mbox{ in }\overline\Omega_i\,,\quad\mbox{for }i\in\mathbb{N}\,.$$ 
Then for every $m\in\mathbb{N}$, we take
$$\mu_m=(1/m,\dots,1/m)\in\Gamma_m^1$$ 
and define the function $$u^\sharp_{p,m}:\overline\Omega^\sharp_m\to[0,+\infty)$$ 
as the $(p,\mu_m)$-mean $u_{p,\mu_m}$ of the functions $u_1,\dots,u_m$, according to (\ref{pconcaveenvelopem}).
\begin{lemma} In the assumptions and notation given above, for $p>0$ the sequence $\{u^\sharp_{p,m}\}_{m=1}^\infty$ is uniformly convergent (up to a subsequence) in $\overline\Omega^\sharp$ to a function $u_p^\sharp\in C(\overline\Omega^\sharp)$ vanishing on $\partial\Omega^\sharp$.
\end{lemma}
{\em Proof. }
First we notice that, by Lemma \ref{proprietaupm}, it follows $u^\sharp_{p,m}\in C(\overline\Omega_m^\sharp)$, $u^\sharp_{p,m}(x)>0$ for $x\in\Omega_m^\sharp$ and
$u^\sharp_{p,m}(x)=0$ for $x\in\partial\Omega_m^\sharp$. 

Since $u^p$ is continuous in the compact set $\overline{\Omega}$, we have that  $u^p$ is uniformly continuous in $\overline\Omega$ and we denote by $\omega_p$ its modulus of continuity.
Obviously $u_i^p$ is also uniformly continuous in $\overline \Omega_i$ with the same modulus of continuity $\omega_p$ for every $i\in\NN$.
Then, as supremal convolution of functions with the same modulus of continuity, also $(u^\sharp_{p,m})^p$ is uniformly continuous in $\overline\Omega_m^\sharp$ with modulus of continuity $\omega_p$ for every $m\in\NN$ (see \cite{st} for instance).
Now let $R=\max\{\mbox{dist}(x,0)\,:\,x\in\overline\Omega\}$ and let $B=B(0,2R)$ the ball centered at the origin with radius $2R$. Then $\rho\overline\Omega\subset B$ for every $\rho\in SO(n)$, so that $\overline\Omega_m^\sharp\subset B$ for every $m\in\NN$ and consequently $\overline\Omega^\sharp\subset B$.
We set $u(x)=0$ for $x\in\overline{B}\setminus\Omega$, $u_i(x)=0$ for $x\in\overline{B}\setminus\Omega_i$ for every $i$ and $u^\sharp_{p,m}(x)=0$ for $x\in\overline{B}\setminus\Omega_m^\sharp$ for every $m$. So extended, $u^p,\,u_i^p$ and $(u^\sharp_{p,m})^p$ obviously remain uniformly continuous, with the same modulus of continuity $\omega_p$, in $\overline B$.
Moreover, since $0\leq u_i(x)\leq M$ for $x\in\rho_i\overline B$, where $M=\max_{\overline\Omega}u$, it follows $0\leq u^\sharp_{p,m}(x)^p\leq M^p$ for $x\in\overline B$ for every $m$. Finally $\{(u^\sharp_{p,m})^p\}_{m=1}^\infty$ is a sequence of equibounded and equicontinuous functions in $\overline B$ and it is possible to extract an uniformly convergent subsequence, say $(u^\sharp_{p,m_k})^p$.  
We set 
\begin{equation}\label{upstardef}
u_p^\sharp(x)=\left(\lim_{k\to+\infty}u^\sharp_{p,m_k}(x)^p\right)^{1/p}\quad x\in\overline B.
\end{equation}
We have just to prove that $u_p^\sharp(x)=0$ for $x\in\partial\Omega^\sharp$.  Since $u_p^\sharp\in C(\overline B)$ (as uniform limit of a sequence of continuous functions), we can just prove that 
$u_p^\sharp$ vanishes  in $B\setminus\overline\Omega^\sharp$. Then consider a point $x\in B\setminus\overline\Omega^\sharp$, that is a point $x$ such that $\mbox{dist}(x,\Omega^\sharp)=d>0$: since $\Omega_m^\sharp$ converges to 
$\Omega^\sharp$ in Hausdorff metric as $m\to\infty$, there exists $M_x$ such that $d(x,\Omega_m^\sharp)>d/2$ for every $m\geq M_x$. Then $x\in B\setminus\Omega_m^\sharp$ and $u^\sharp_{p,m}(x)=0$ for $m\geq M_x$, whence $u_p^\sharp(x)=0$. 
\qed
\medskip

The previous lemma contains the definition  (\ref{upstardef}) of $u_p^\sharp$, which may look quite involved, as already said in the Introduction. 
To give a geometric insight, let me say that it is somewhat reminiscent of a rearrangement technique introduced by Tso in \cite{Tso} to treat the case of Monge-Amp\`ere equation, where every sublevel set of a convex function is substituted by a ball with the same mean width; here, instead, the level sets of $u_p^\sharp$ are not necessarily balls, apart from $\Omega^\sharp=\{u_p^\sharp\geq 0\}$, and their mean width is in general greater than the mean width of the corresponding level sets of $u$, apart again from the ground domain (indeed, we are not acting separately on each level sets, but globally on the function).
Notice also that (\ref{upstardef}) may be considered not completely satisfying as a definition of a rearrangement, since it seems to depend on the chosen subsequence $u_{p,m_k}^\sharp$. 
Nevertheless 
it suffices to prove 
a priori estimates similar to (\ref{comparisonTalenti})-(\ref{comparisonTalenti1}) of the solution ${u}$ of (\ref{pb1}) 
in terms of the solution $v$ in $\Omega^\sharp$  when $F$ is a {\em rotationally invariant} operator,
i.e. when
\begin{equation}\label{Frotinvariant}
F(\rho x, u, \rho\theta,\rho A\rho^T)=F(x,u,\theta,A)
\end{equation}
for every $(x,u,\theta,A)\in\RR^n\times\RR\times\RR^n\times S_n$ and every $\rho\in SO(n)$. 

Examples of rotationally invariant operators are the Laplacian, the $q$-Laplacian, the mean curvature operator, the Hessian operators, etc.
Moreover notice that $F$ is rotationally invariant when it depends on $x,\,\theta$ and $A$ only in terms of $|x|$, $|\theta|$ and the eigenvalues of $A$, respectively. 

\begin{remark}
If $F$ is rotationally invariant and $u$ solves (\ref{pb1}) in $\Omega$, then {\em $u_\rho$ (defined in (\ref{urho})) solves (\ref{pb1}) in $\rho\Omega$.}
This is the reason for we consider rotationally invariant operators.
\end{remark}

In view of  (\ref{Gpno0}), given an operator $F$, a real number $p>0$ and a vector $\theta\in\RR^n$, we set
\begin{equation}\label{Gptheta}
G_{p}^{(\theta)}(x,t,A)= F(x,t^{\frac{1}{p}},t^{\frac{1}{p}-1}\theta,
t^{\frac{1}{p}-3}A)\qquad(x,t,A)\in\RR^n\times[0,\infty)\times S_n\,.
\end{equation}
 
\begin{lemma}\label{lemmarearr}
Let $\Omega$ be a bounded open convex set in $\RR^n$ and $u$ a classical solution of (\ref{pb1}) in $\Omega$,
where $F$ is a rotationally invariant proper elliptic operator and $u$ satisfies assumption (\ref{liminf}) for every $x\in\partial\Omega$.

Let $p\in(0,1)$ and assume that  
\begin{equation}\label{AssGp0conv}
\mbox{the set }\,\{(x,t,A)\in[0,\infty)\times S_n\,:\,G_p^{(\theta)}(x,t,A)\geq 0\}\,\,\mbox{ is convex}
\end{equation}
for every fixed $\theta\in\RR^n$.

Then $u_p^\sharp$ is a viscosity subsolution of problem (\ref{pb1}) in $\Omega^\sharp$.
\end{lemma}
{\em Proof. }
Set $\mu=(1/m,\dots,1/m)$, $F_\mu=F$ and $F_i=F$ for $i=1,\dots, m$. 
Then assumption (\ref{AssGp0conv}) implies that ($WC_{\mu,p}$) holds for the operators $F_\mu$ and $F_1,\dots, F_m$.
Hence $u^\sharp_{p,m}$ is a viscosity subsolution in $\Omega_m$ for every $m\in\NN$ by Lemma \ref{lemmapleq02}.
The conclusion follows thanks to the stability of viscosity subsolutions with respect to uniform convergence.
\qed\medskip

\begin{remark}\label{hpGpquasiconcave}
Notice that assumption (\ref{AssGp0conv}) is satisfied if the function $G_p^{(\theta)}$ is quasi-concave for every $\theta\in\RR^n$, hence if it is $q$-concave for some $q\in\RR$.
\end{remark}

\begin{theorem}\label{thmrearr}
In the same assumptions of the previous lemma, if $F$ satisfies the comparison principle $(CP)$, then (\ref{comparisonSalani}) and (\ref{lpustargequ}) hold, where 
${u}$ and ${v}$ are the solutions of problem (\ref{pb1}) and (\ref{pbsharp}), respectively.
\end{theorem}
{\em Proof. }
First notice that (\ref{comparisonSalani}) follows from Lemma \ref{lemmarearr} using the comparison principle.

Then, (\ref{comparisonSalani}) yields 
\begin{equation}\label{Lpcomparison1}
\|u_p^\sharp\|_{L^q(\Omega^\sharp)}\leq\|{v}\|_{L^q(\Omega^\sharp)}\quad\mbox{for every }q\in(0,\infty]\,.
\end{equation}

Next we prove that for every $q\in(0,\infty]$
\begin{equation}\label{Lpcomparisonm}
\|{u}\|_{L^q(\Omega)}\leq
\|u^\sharp_{p,m}\|_{L^q(\Omega_m)}\quad\mbox{for every }m\in\NN\,.
\end{equation}
The case $q=+\infty$ is almost trivial, since we obviously have $M=\max_{\overline\Omega}{u}=\max_{\overline\Omega_i}{u_i}$ for every $i$, then
$\max_{\overline\Omega^\sharp_m} u^\sharp_{p,m}=M$ for every $m$.
Now let $q\in(0,\infty)$. By the layer cake formula, it holds
\begin{equation}\label{layer1}
\|{u}\|_{L^q(\Omega)}=\int_0^M|\Omega(t)|\,dt\,,
\end{equation}
where
$$
\Omega(t)=\{x\in\Omega\,:\,{u}(x)\geq t\}\,,
$$
and
\begin{equation}\label{layer2}
\|u^\sharp_{p,m}\|_{L^q(\Omega_m)}=\int_0^M|\Omega_m(t)|\,dt\,,
\end{equation}
where
$$
\Omega_m(t)=\{x\in\Omega_m\,:\,u^\sharp_{p,m}(x)\geq t\}\,.
$$
On the other hand, the definition of $u^\sharp_{p,m}$ implies
$$
\Omega_m(t)\supseteq\frac{1}{m}\sum_{i=1}^m\rho_i\Omega(t)\,,
$$
then 
\begin{equation}\label{notequidistributed}
|\Omega_m(t)|\geq|\Omega(t)|\qquad\mbox{for every }t\in[0,M]
\end{equation}
by the Brunn-Minkowski inequality (\ref{BMineq}). In view of (\ref{layer1}) and (\ref{layer2}), (\ref{notequidistributed}) yields (\ref{Lpcomparisonm}).

Finally, since 
$$
\|u_p^\sharp\|_{L^q(\Omega^\sharp)}=\lim_{m\to\infty}\|u^\sharp_{p,m}\|_{L^q(\Omega_m)}
$$
by uniform convergence, then by coupling (\ref{Lpcomparison1}) and
(\ref{Lpcomparisonm}) we get (\ref{lpustargequ})  and the proof is complete.
\qed\medskip

As a final remark on $u_p^\sharp$, let me notice that it is not equidistributed with ${u}$ (in contrast with Schwarz symmetrization), but in fact each one of its level sets has greater measure than the corresponding level set of ${u}$ by (\ref{notequidistributed}).

\section{Examples}
We first give examples of applications of Theorem \ref{mainthm}, then we consider the mean width rearrangement and give examples of applications of Theorem \ref{thmrearr}. Throughout this section, $f$ and $f_i$, $i=0,1,\mu$, are non-negative smooth functions.

\subsection{Applications of Theorem \ref{mainthm}.}
The first example is the Laplacian. Let:
$$
F_i(x,u,Du,D^2u)=\Delta u + f_i(x,u,Du)\qquad i=0,1,\mu
$$
and
\begin{equation}\label{gip}
g_{i,p}^{(\theta)}(x,t)=\left\{\begin{array}{lll}t^{3-\frac{1}{p}}f(x,t^{\frac{1}{p}},t^{\frac{1}{p} - 1}\theta)\quad&p>0\\
&&\quad i=0,1,\mu\,.\\
e^{-t}f(x,e^t,e^t\theta)\quad&p=0
\end{array}\right.
\end{equation}
In this case condition $(WA_{\mu,p})$ is satisfied if
\begin{equation}\label{hp4laplacian}
g_{\mu,p}^{(\theta)}((1-\mu)x_0+\mu x_1, (1-\lambda)t_0+\lambda t_1)\geq(1-\lambda)g_{0,p}^{(\theta)}(x_0,t_0)+\lambda g_{1,p}^{(\theta)}(x_1,t_1)
\end{equation}
for every $x_0\in\Omega_0$, $x_1\in\Omega_1$, $t_0,t_1\geq 0$ and for every fixed $\theta\in\rn$.

For instance, let $u_0$ and $u_1$ be the solutions of the following problems
$$
\left\{\begin{array}{ll}
\Delta u_0+f_0(x)=0\quad&\mbox{in }Q=[-1,1]\times[-1,1]\\
\\
u_0=0\quad&\mbox{on }\partial Q
\end{array}\right.
$$
and
$$
\left\{\begin{array}{ll}
\Delta u_1+f_1(x)=0\quad&\mbox{in }B(0,1)\\
\\
u_1=0\quad&\mbox{when }|x|=1,
\end{array}\right.
$$
respectively.

\begin{figure}[htbp]
\begin{center}
\includegraphics[scale=.4]{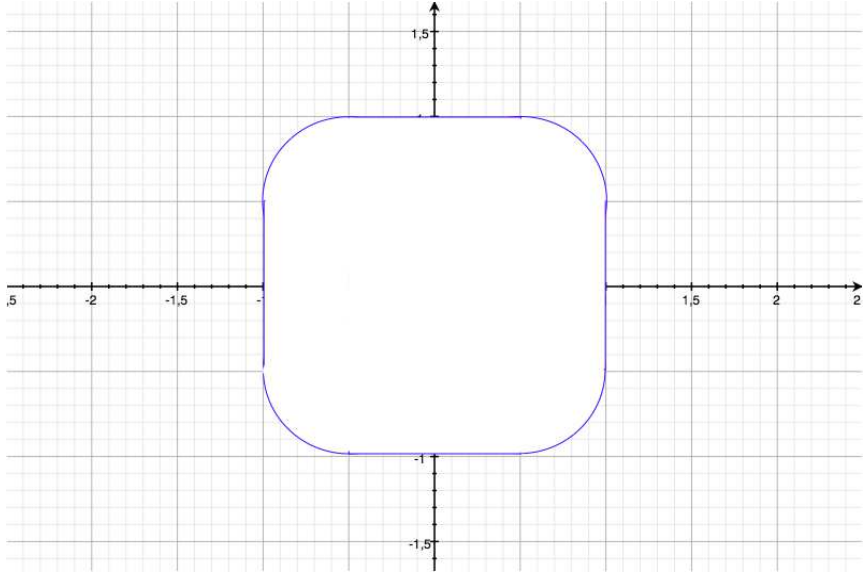}
\caption{The Minkowski combination of a square and a circle: $\Omega=\frac{1}{2} Q+\frac12 B(0,1)$.}
\label{fig1}
\end{center}
\end{figure}

Then take $\mu=1/2$ and set
$$
\Omega=\frac{1}{2} Q+\frac12B(0,1)\,,
$$
see Figure \ref{fig1}.
Now let $u_\mu$ be the solution of
$$
\left\{\begin{array}{ll}
\Delta u_\mu+f_\mu(x)=0\quad&\mbox{in }\Omega\\
\\
u_\mu=0\quad&\mbox{on }\partial\Omega.
\end{array}\right.
$$
Then (\ref{hp4laplacian}) for p=$1/3$ reads
\begin{equation}\label{nuovahp}
f_\mu\big(\frac{x_0+x_1}{2}\big)\geq\frac12f_0(x_0)+\frac12f_1(x_1)
\end{equation}
(please, compare with (\ref{eq:4.2})) and
Theorem \ref{mainthm} tells that we can estimate $u_\mu$ in terms of $u_0$ and $u_1$; precisely it holds
$$
u_\mu\big(\frac{x_0+x_1}{2}\big)\geq\left[\frac12\sqrt[3]{u_0(x_0)}+\frac12\sqrt[3]{u_1(x_1)}\right]^3\quad\mbox{for every }x_0\in Q,\,x_1\in B(0,1)
$$
and Corollary \ref{maincor} yields
$$
\|u_\mu\|_{L^r(\Omega_\mu)}\geq M_q(\|u_0\|_{L^r(Q)},\|u_1\|_{L^r(B(0,1))};\mu)
$$
for every $r\in(0,+\infty]$, where 
$$
q=\left\{\begin{array}{ll}\frac{r}{n+3r}\,,\qquad&r\in(0,+\infty)\\
\\
1/3\,,\qquad&r=+\infty
\end{array}\right.\,.
$$
Notice in particular that, if $$f_0=f_1=f_\mu=f:\rn\to[0,+\infty)\,,$$ condition (\ref{nuovahp}) simply means {\em $f$ is concave}.
More generally, in this particular case, we can write the following result.
\begin{proposition}
Let $f$ be a smooth nonnegative function defined in $\rn$. Let $\mu\in(0,1)$ and $\Omega_0$ and $\Omega_1$ be convex subsets of $\rn$ and
denote by $u_0$, $u_1$ and $u_\mu$ the solutions of
$$
\left\{\begin{array}{ll}
\Delta u_i+f(x)=0\quad&\mbox{in }\Omega_i\\
\\
u_i=0\quad&\mbox{on }\partial\Omega_i
\end{array}\right.
$$
for $i=0,1,\mu$ respectively, where $\Omega_\mu=(1-\mu)\Omega_0+\mu\Omega_1$, as usual.

Assume $f$ is $\beta$-concave for some $\beta\geq 1$, that is $f^\beta$ is concave.

Then (\ref{eqmainthm}) holds with 
$$
p=\frac{\beta}{1+2\beta}
$$
and consequently (\ref{eqLpnorm}) holds with 
$$
q=\left\{\begin{array}{ll}\frac{\beta r}{n\beta+r(1+2\beta)}\,,\qquad&\mbox{for }r\in(0,\infty)\\
\\
\frac{\beta}{1+2\beta}\,,\qquad&\mbox{for }r=+\infty\,.
\end{array}\right.
$$
In case $f$ is a positive constant ($\beta=+\infty$), the same conclusions follow with $p=1/2$ and 
$$q=\left\{\begin{array}{ll}r/(n+2r)\,,\quad&\mbox{for }r\in(0,\infty)\\
 1/2\,,\quad&\mbox{for }r=+\infty\,.
 \end{array}\right.
 $$
\end{proposition}
{\em Proof. }
The proof is a direct application of Theorem \ref{mainthm2}, in view of Lemma A.1 of \cite{IS}.
\qed\medskip

Notice that the assumptions of the above proposition imply that the involved solutions $u_0$, $u_1$ and $u_\mu$ are all $p$-concave in their own domains,
see \cite{BS, Kenni, Kore2}.
\medskip

The example of the Laplace operator can be generalized by considering the $q$-laplacian for $q>1$:
$$
F(x,u,Du,D^2u)=\Delta_q u + f(x,u,Du)\,,
$$
where $\Delta_q u=\mbox{div}(|Du|^{q-2}Du)$, as usual.
In this case we set
\begin{equation}\label{conditionqLaplacian}
g_{i,p}^{(\theta)}(x,t)=\left\{\begin{array}{lll}
t^{q+1-\frac{q-1}{p}}f_i(x,t^{\frac{1}{p}},t^{\frac{1}{p}-1}\theta)\quad&\mbox{if }p\neq 0\\
&&\quad i=0,1,\mu\,.\\
e^{t(1-q)}f_i(x,e^t,e^t\theta)\quad&\mbox{if }p=0
\end{array}\right.
\end{equation}
and we can apply Theorem \ref{mainthm2} again if (\ref{hp4laplacian}) holds.
\medskip

Another generalization of the Laplacian is
the \emph{Finsler Laplacian} $\Delta_Hu$, which for a regular function $u$ is defined as follows
$$
\Delta_H u=\mbox{div}\large(H(Du)\nabla_\xi H(Du)\large)\,,
$$
where $H(\xi)$ is a given norm in $\RR^n$, that is a nonnegative centrally symmetric $1$-homogeneous convex function (or, if you prefer, the support function of a centrally symmetric convex body), and $\nabla_\xi$ denotes the gradient with respect to the variable $\xi\in\RR^n$. For more detail, please refer for instance to \cite{FK} and references therein.
Our results can be applied to the operator
$$
F(x,u,Du,D^2u)=\Delta_H u + f(x,u,Du)
$$
exactly in the same assumptions as for the Laplacian, that is when (\ref{hp4laplacian}) holds, where $g_{i,p}^{(\theta)}$ is given by (\ref{gip}).


The Laplacian and the $q$-Laplacian are however classical matters of investigation and the results above
stated could be considered to some extent not completely new. Indeed results similar to the ones we could obtain by applying Theorem \ref{mainthm2} to some suitable particular situation, have been already used to prove Brunn-Minkowski type inequalities for some variational functionals: see for instance \cite{Co, CoCuSa} for the Laplacian and $q$-Laplacian, while in the case on the Finsler Laplacian related results can be found in \cite{WX}. 
\medskip

Completely new applications are instead obtained when considering for instance
Dirichlet problems for
the Pucci's Extremal Operator $\mathcal{M}^-_{\lambda, \Lambda}u$.

The {\em Pucci's Extremal Operators} were introduced by C. Pucci
in \cite{Pucci} and they are perturbations of the usual Laplacian.
Precisely, given two numbers $0<\lambda \leq \Lambda$ and a real symmetric $n
\times n$ matrix $M$, whose eigenvalues are $e_i=e_i(M)$, for
$i=1,...,n$, the Pucci's extremal operators are
\begin{equation}\label{Pucci+} {\mathcal
M}^+_{\lambda,\Lambda}(M)=\Lambda
\sum_{e_i>0}e_i+\lambda\sum_{e_i<0}e_i
\end{equation}
and
\begin{equation} \label{Pucci-} {\mathcal
M}^-_{\lambda,\Lambda}(M)=\lambda
\sum_{e_i>0}e_i+\Lambda\sum_{e_i<0}e_i.\end{equation}
We recall that ${\mathcal
M}^+_{\lambda,\Lambda}$ and ${\mathcal
M}^-_{\lambda,\Lambda}$ are uniformly elliptic and positively homogeneous of degree $1$; moreover ${\mathcal
M}^+_{\lambda,\Lambda}$ is convex, while ${\mathcal
M}^-_{\lambda,\Lambda}$ is concave over $S_n$ (see \cite{caffacabre} for instance). 

Again Theorem \ref{mainthm2} can be applied to
$$F(x,u,Du,D^2u)=\mathcal{M}^-_{\lambda, \Lambda}(D^2u) + f(x,u,D u)\,
$$
in the same assumptions as for the Laplacian and the Finsler laplacian, that is when (\ref{hp4laplacian}) holds, where $g_{i,p}^{(\theta)}$ is given by (\ref{gip}).

In fact the same conclusion holds for every elliptic equation of the type
$$
F(D^2u)+ f(x,u,Du)=0\,,
$$
where $F:S_n\to\RR$ is concave and positively $1$-homogeneous, as it is easily seen.

\subsection{Applications of Theorem \ref{thmrearr}.}

Next we discuss applications of the rearrangement technique introduced in Section \ref{rearrsect} to the examples given in the previous subsection.
Also in light of the above discussion, it is easily seen that Theorem \ref{thmrearr} can be applied to the Laplacian  and in particular to the problem
$$
\left\{\begin{array}{ll}
\Delta {u} + f(|x|,u,|Du|)=0\quad&\mbox{in }\Omega\,,\\
{u}=0\quad&\mbox{on }\partial\Omega\,,
\end{array}\right.
$$
when the function
$$
g_p(s,t)=t^{3-\frac{1}{p}}f(s,t^{\frac{1}{p}},t^{\frac{1}{p} - 1}r)
$$
is not increasing with respect to $s\in[0,+\infty)$ and concave with respect to $(s,t)\in[0,\infty)^2$, for every $r\geq 0$, for some $p\in(0,1)$. In this case we can compare the mean width rearrangement of ${u}$ with the solution ${v}$ of the same problem in the ball $\Omega^\sharp$ wit the same mean width of $\Omega$ and finally get (\ref{lpustargequ}).

For instance, when $f=1$, we can express the result so obtained in a striking way in terms of torsional rigidity: {\em among convex sets of given mean width, the torsional rigidity is maximized by the ball}. When we denote by $\tau(A)$ the torsional rigidity of the set $A$, we can translate the latter sentence in the following Urysohn's type inequality (see \cite{GS}):
$$
\tau(\Omega)\leq\tau(\Omega^\sharp)\quad\mbox{ for every convex set }\Omega\,.
$$
On the other hand this is weaker than the well known inequality (see \cite{PS}) $\tau(\Omega)\leq\tau(\Omega^\star)\,,$ since $\tau$ is increasing with respect to inclusion and 
$$\Omega^\star\subseteq\Omega^\sharp$$ 
by the classical Urysohn's inequality between mean width and volume. 
The latter in fact implies that, in general, we cannot expect to find interesting new inequalities by applying Theorem \ref{thmrearr} to equations involving the Laplacian or some other divergence type operator, which Schwarz symmetrization fits better.
Similar considerations obviously hold for the $p$-Laplacian.
In the case of Finsler Laplacian instead, Theorem \ref{thmrearr} neither can be  applied due to the fact that $\Delta_H$ is not invariant with respect to rotation.
\medskip

When instead considering the Pucci's extremal operators, Schwarz symmetrization has not been successfully applied until now, at least to my knowledge. In this case Theorem 
\ref{thmrearr} can be applied to the operator $$F(x,u,Du,D^2u)=\mathcal{M}^-_{\lambda, \Lambda}(D^2u) + f(|x|,u,|D u|)$$
precisely under the same assumptions as for the Laplacian, 
and in this case it  yields completely new results. 

As a paradigmatic explicit example, let me write the following proposition.
\begin{proposition}
Let $\Omega$ be an open bounded convex set in $\rn$ and $\Omega^\sharp$ a ball with the same mean width of $\Omega$.
Let $u$ and $v$ solve the following problems
$$
\left\{\begin{array}{ll}
\mathcal{M}^-_{\lambda, \Lambda}(D^2u)+1=0\quad&\mbox{ in }\Omega\,,\\
u=0\quad&\mbox{ on }\partial\Omega\,,
\end{array}\right.
$$
and
$$
\left\{\begin{array}{ll}
\mathcal{M}^-_{\lambda, \Lambda}(D^2v)+1=0\quad&\mbox{ in }\Omega^\sharp\,,\\
v=0\quad&\mbox{ on }\partial\Omega^\sharp\,,
\end{array}\right.
$$
respectively.

Then 
$$
\|u\|_{L^p(\Omega)}\leq\|v\|_{L^p(\Omega^\star)}
$$
for every $p>0$, including $p=+\infty$.
\end{proposition}

Similarly to what observed at the end of the previous subsection, the results obtained for the Pucci operator $\mathcal{M}^-_{\lambda, \Lambda}$ hold in fact for every elliptic equation of the type
$$
F(D^2u)+ f(|x|,u,|Du|)=0\,,
$$
where $F:S_n\to\RR$ is concave, rotationally invariant and positively $1$-homogeneous, as it is easily seen.


\end{document}